\documentclass[11pt,reqno,twoside]{article}

\usepackage{paper_style}

\usepackage{mathtools}
\usepackage{graphicx}
\usepackage{multirow}
\usepackage{amsmath,amssymb,amsfonts}
\usepackage{mathrsfs}
\usepackage{url}
\usepackage{xcolor}
\usepackage{textcomp}
\usepackage{manyfoot}
\usepackage{booktabs}
\usepackage{xspace}
\usepackage{rotating}
\usepackage{stmaryrd}
\usepackage{scalefnt}
\usepackage{upgreek}

\definecolor{keywordcolor}{rgb}{0.7, 0.1, 0.1}   
\definecolor{tacticcolor}{rgb}{0.0, 0.1, 0.3}    
\definecolor{commentcolor}{rgb}{0.4, 0.4, 0.4}   
\definecolor{stringcolor}{rgb}{0.5, 0.3, 0.2}    
\definecolor{symbolcolor}{rgb}{0.1, 0.2, 0.7}    
\definecolor{sortcolor}{rgb}{0.1, 0.5, 0.1}      
\definecolor{attributecolor}{rgb}{0.7, 0.1, 0.1} 
\definecolor{errorcolor}{rgb}{1, 0, 0}           

\usepackage{listings}

\newcommand{\lean}[1]{\lstinline[language=lean, mathescape=true]{#1}}

\newcommand{\mathlib}{\textsf{mathlib}\xspace}

\newcommand{\Fix}{\mathrm{Fix}\,}

\newcommand{\Hh}{\mathcal{H}}
\newcommand{\ip}[2]{\langle #1, #2 \rangle}
\newcommand{\email}[1]{\protect\href{mailto:#1}{#1}}

\newcommand{\KM}{Krasnosel'ski\u{\i}--Mann\xspace}
\newcommand{\Fejer}{Fej\'{e}r\xspace}

\begin{document}
\title{Formalization of Two Fixed-Point Algorithms in Hilbert Spaces} 

\author{Yifan Bai\thanks{School of Mathematical Sciences, Shanghai Jiao Tong University, Shanghai China (\email{yifanbai@sjtu.edu.cn}).}
\and Yantao Li\thanks{Zhiyuan College, Shanghai Jiao Tong University, Shanghai China (\email{alexandero@sjtu.edu.cn}).}
\and Jian Yu\thanks{School of Mathematics, Jilin University, Jilin China (\email{yujian6521@mails.jlu.edu.cn}).}
\and Jingwei Liang\thanks{School of Mathematical Sciences \& Institute of Natural Sciences, Shanghai Jiao Tong University, Shanghai China (\email{jingwei.liang@sjtu.edu.cn}).}}
\date{}
\maketitle



\begin{abstract}
Iterative algorithms are fundamental tools for approximating fixed-points of nonexpansive operators in real Hilbert spaces. Among them, Krasnosel'ski\u{\i}--Mann iteration and Halpern iteration are two widely used schemes. In this work, we formalize the convergence of these two fixed-point algorithms in the interactive theorem prover Lean4 based on type dependent theory. To this end, weak convergence and topological properties in the infinite-dimensional real Hilbert space are formalized. Definition and properties of nonexpansive operators are also provided. As a useful tool in convex analysis, we then formalize the Fej\'{e}r monotone sequence. Building on these foundations, we verify the convergence of both the iteration schemes. Our formalization provides reusable components for machine-checked convergence analysis of fixed-point iterations and theories of convex analysis in real Hilbert spaces. Our code is available at \url{https://github.com/TTony2019/fixed-point-iterations-in-lean}.
\end{abstract}

\paragraph{Key words} Formalization, Lean, Fixed-point Iteration, Halpern Iteration, \KM Iteration

\section{Introduction}
\label{sec:introduction}

Fixed-point iterative schemes, which aims to find a fixed-point of some given operator, plays an important role in modern computational mathematics including numerical analysis, scientific computing and optimization, {\it etc}. 
In this paper, we consider the fixed-point problem of the following form
\begin{equation}
\label{eq:main_prob}
    \text{Find}~~ x \in \Hh \quad{\rm such~that}\quad x \in \Fix T \cap D,
\end{equation}
where $\Hh$ is a real Hilbert space, $D\subset \Hh$ is a nonempty closed convex set, $T : \Hh \to \Hh$ is a nonexpansive operator with $T(D) \subseteq D$ and $\Fix T = \{ x \in \Hh \mid Tx = x \}$ denotes the set of fixed-points of $T$.

\subsection{Iterative schemes}
Over the past decades, various iterative schemes have been proposed to solve \eqref{eq:main_prob}. 
\textbf{Picard iteration}~\cite{picard1890memoire} is a classical approach, which repeatedly applies the operator $T$:
\begin{equation}\label{eq:picard}
    \forall n \in \mathbb{N}, \qquad x_{n+1} = T x_n.
\end{equation}
The convergence of \eqref{eq:picard} can be guaranteed when $T$ is contractive, by virtue of the Banach contractive mapping principle. However, as mentioned in~\cite{bauschke2011convex}, when $T$ is simply nonexpansive (\cref{def:nonexp}), e.g. $T=\Id$ being the identity operator of $\Hh$, Picard iteration fails to converge.

In the literature, several remedies are proposed \cite{bauschke2011convex}. 
\textbf{Halpern iteration}~\cite{halpern1967fixed} chooses one fixed anchor point $u \in D$, and computes the new iterate by a convex combination of the anchor and the image of the current iterate:
\begin{equation}\label{eq:Halpern}
    \forall n \in \mathbb{N}, \qquad x_{n+1} = \alpha_n u + (1-\alpha_n)Tx_{n},
\end{equation}
where $\{ \alpha_n\} \subseteq [0,1]$. 
Under certain conditions on the weight $\{\alpha_n\}$, Halpern iteration converges {\it strongly} to $P_{\Fix T \cap D}(u)$~\cite{halpern1967fixed}, 
where $P_{\Fix T \cap D}(\cdot)$ denotes the projection operator of the set $\Fix T \cap D$. 
Another well-known approach is the \textbf{Krasnosel'ski\u{\i}--Mann (KM)} iteration~\cite{krasnosel1955two,mann1953mean}, which relaxes Picard iteration by averaging the current iterate and its image under $T$,
\begin{equation}\label{eq:KM}
    \forall n \in \mathbb{N}, \qquad x_{n+1} = (1-\alpha_n) x_n +\alpha_nTx_{n},
\end{equation}
where $\{ \alpha_n\} \subseteq [0,1]$. 
In the realm of convex optimization and set-valued analysis, many numerical schemes can be formulated as instances of KM iteration, including a wide range of operator splitting method, see~\cite{bauschke2011convex} for examples. Early investigations of convergence of KM iteration focus on constant weights in Banach space, see~\cite{krasnosel1955two, schaefer1957uber}. Weak convergence of the general KM iteration \eqref{eq:KM} is considered in ~\cite{reich1979weak}. For more fixed-point type iterative schemes, we refer to~\cite{berinde2007iterative} and the references therein.

\subsection{Motivation and main contributions}

The formalization of mathematics, which plays a crucial rule in AI4Math, has received increasing amount of attention. 
Available formalization languages include Rocq (Coq)~\cite{coq1996coq},  Isabelle~\cite{nipkow2002isabelle}, Lean4 \cite{de2015lean}, {\it etc}. 
In particular, Lean4 is the main formalization language adopted by the mathematics community due to its rapidly growing community library \mathlib~\cite{The_mathlib_Community_2020}. Leveraging this powerful library, significant progress has been made in formalizing convex optimization theory and numerical algorithms in Lean4. For example,  in nonsmooth optimization and first-order optimization methods, \cite{li2025formalization2,li2025formalization1,li2025formalization3} formalize the first-order optimality, convergence analysis of first-order method and convergence rates, respectively.

\begin{figure}[h]
    \centering
    \includegraphics[width=0.96\linewidth]{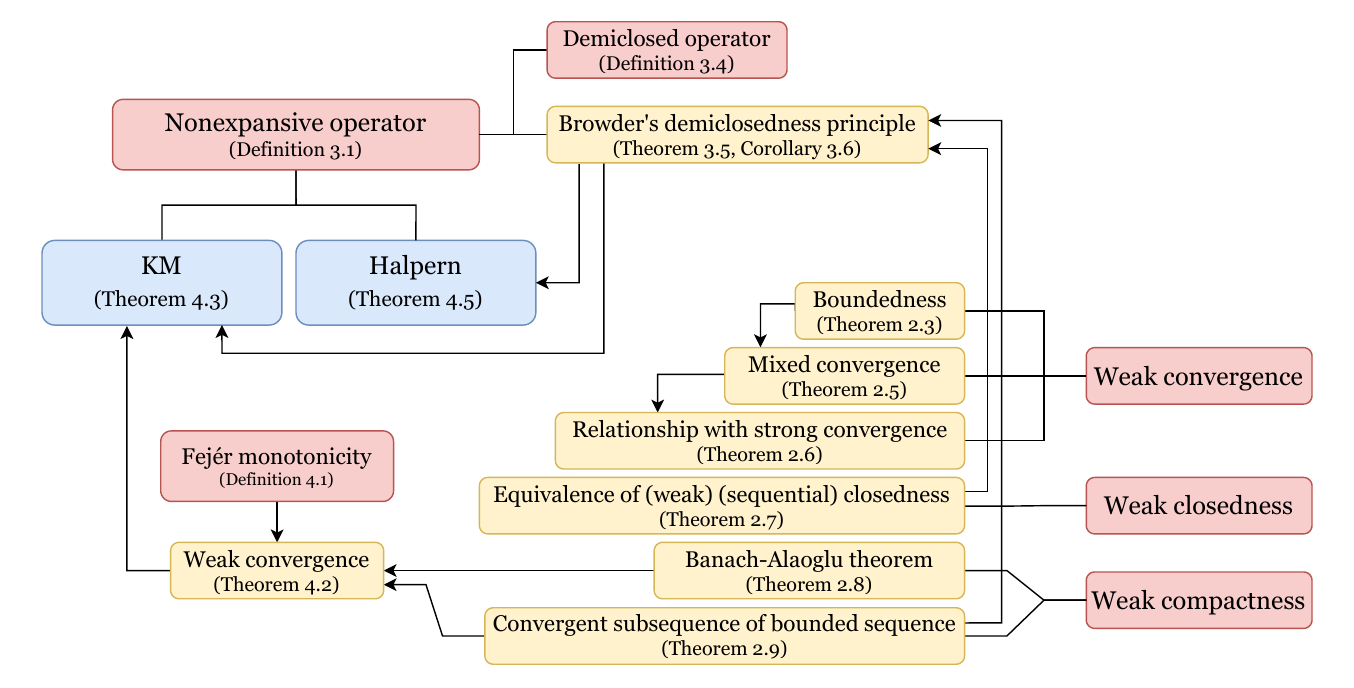}
    \caption{Dependency graph of part of the formalized theorems.}
    \label{fig:paper_struct}
\end{figure}

The formalization of iterative schemes has been considered in the literature, such as: in Rocq, the Jacobi method for solving linear systems is verified in~\cite{tekriwal2023verified} and convergence of stationary iterative methods is formalized in~\cite{tekriwal2024formalization}. While in Isabelle, a framework for verifying numerical algorithms (e.g., bisection method and fixed-point method for contractive real functions) is developed in~\cite{bryant2025verifying}.
Given the fundamental role of fixed-point iteration, 
in this paper we use Lean4 to formalize the convergence of Halpern and \KM fixed-point iterations provided in~\cite{bauschke2011convex}. 
Below we list the main challenges to achieve the goal
\begin{itemize}
    \item Limited result on \lean{WeakSpace} available in \mathlib (see \cref{sec:pre}). Also, there exists a mismatch between the net-based convergence in \cite{bauschke2011convex} and the filter-based convergence framework used in \mathlib.

    \item For nonexpansive operator and related result, they are currently unavailable in \mathlib. 
\end{itemize}
To achieve the formalization goal, the necessary results and their dependency relation are provided in \cref{fig:paper_struct}: the fundamentals are the weak topology in infinite-dimensional real Hilbert space (see \cref{sec:Weak}), related definition and properties associated to nonexpansive operator (see \cref{sec:nonexp}), and lastly the convergence of \KM and Halpern iterations (see \cref{sec:convergence}).

\paragraph{Main contributions}
By formalizing the convergence of two widely adopted fixed-point iterative schemes, in this paper our contribution consists of the following aspects
\begin{enumerate}
    \item We formalize the weak convergence characterized by the inner product in real Hilbert spaces and demonstrate its equivalence to the filter-based definition of convergence in \mathlib. While several related results in \mathlib are formalized in the weak-\(*\) topology, our formalization works directly in weak topology. This provides a more convenient and reusable interface for subsequent formalizations in functional analysis and optimization.

    \item We provide a formalization of the nonexpansive operator and a collection of associated key properties. In particular, we formalize the relationship between nonexpansive operators and demiclosed operators (e.g., Browder's demiclosedness principle). This lays the groundwork for formally verifying a broader class of optimization and operator splitting algorithms.
    
    \item We formalize several key properties of Fej\'{e}r monotonicity related to the nonexpansive operator, which is widely used for the convergence analysis of the fixed-point iteration. Our formalization includes lemmas that connect Fej\'{e}r monotonicity with boundedness and convergence properties of sequences with respect to closed convex sets, making it a reusable tool for later algorithmic analysis.
    
    \item Building on the above infrastructure, we formalize the convergence analysis of two classical fixed-point iterations: the Halpern iteration and the Krasnosel'ski\u{\i}--Mann (KM) iteration. These results demonstrate how the formalized weak topology toolkit, nonexpansive operator theory, and Fej\'{e}r monotonicity can be combined to verify nontrivial convergence theorems in Lean4.
\end{enumerate}

\paragraph{Paper organization}
The rest of the paper is organized as follows. A brief review of definitions of weak topology including \lean{WeakDual} and \lean{WeakSpace} is provided in \cref{sec:pre}. 
Formalization of the properties on the weak topology of Hilbert spaces are presented in \cref{sec:weak}. We formalize the definition of nonexpansive operators in \cref{sec:nonexp} and Fej\'{e}r monotone sequence in \cref{sec:fejer}. The convergence analysis of Halpern iteration and KM iteration are formalized in \cref{sec:halpern} and \cref{sec:KM} respectively. 
Finally, we conclude this paper and remark for future research directions in \cref{sec:con}.

\section{Weak topology of Hilbert spaces}
\label{sec:Weak}
In this section, we develop the basic weak-topological framework for Hilbert spaces. A main challenge in our formalization is the current lack of supporting results about \lean{WeakSpace} in \mathlib. To highlight this gap, we first summarize the relevant material already available in \mathlib, and then present the additional definitions and lemmas established in our development.

\subsection{Existing results in \mathlib}
\label{sec:pre}

We first introduce some preliminaries of weak topology over a general Hilbert space, and the mathematical background of the corresponding formalization in \mathlib.

\begin{definition}[Weak topology]
\label{def_weak}
    For a commutative semiring $k$ and a $k$-module $E$,  the weak-\(*\) topology on the (topological) dual $E^*$ is defined as the coarsest topology such that for all $x \in E$, the map
\begin{equation}
  f_x :\; E^* \to k ,\qquad \phi \mapsto \phi(x)
\end{equation}
is continuous. On the other hand, the weak topology on $E$ is the coarsest topology such that for all $\phi \in E^*$, the map
\begin{equation}
  f_\phi :\; E \to k ,\qquad x \mapsto \phi(x)
\end{equation}
is continuous.
\end{definition}

In \mathlib, the weak topologies between $k$-modules $E$ and $F$ are introduced in a more general setting via bilinear pairings. Let
\begin{equation}
\label{eq:def_B}
    B : E\times F \to k
\end{equation}
be a bilinear mapping. For any $y \in F$, the mapping $x \mapsto B(x,y)$ is a $k$-linear mapping from $E$ to $k$. The \emph{weak topology on $E$ induced by $B$} is the coarsest topology such that $x \mapsto B(x,y)$ is continuous for all $y \in F$, which is formalized by \lean{WeakBilin}:
\begin{lstlisting}[frame=single]
def WeakBilin [CommSemiring k] [AddCommMonoid E] [Module k E] [AddCommMonoid F] [Module k F] (_ : E →ₗ[k] F →ₗ[k] k) := E
\end{lstlisting}

The type of the bilinear pairing $B$ (\cref{eq:def_B}) is defined in a currying form \lean{E →ₗ[k] F →ₗ[k] k}, and the weak topological space \lean{WeakBilin} over \lean{E} is a type synonym: it has the same underlying type as the original space \lean{E}, but is intended to carry a different \lean{TopologicalSpace} instance—namely, the weak topology induced by a bilinear pairing \lean{B}.

To recover \cref{def_weak} from this general definition, one needs to apply \lean{WeakBilin} to the canonical pairing between $E$ and its topological dual $E^*$ (in \mathlib, \lean{E →L[k] k}).
Consider the topological dual pairing
\begin{equation}
    B : E^* \times E \to k, \qquad \ip{\phi}{x} \mapsto \phi(x).
\end{equation}
In \mathlib, this canonical pairing is formalized as
\begin{lstlisting}[frame=single]
def topDualPairing : (E →L[k] k) →ₗ[k] E →ₗ[k] k := ContinuousLinearMap.coeLM k
\end{lstlisting}

Based on the canonical pairing, the weak-\(*\) topology on $E$ is equivalent to \lean{WeakBilin (topDualPairing k E)}:
\begin{lstlisting}[frame=single]
def WeakDual (k E : Type*) [CommSemiring k] [TopologicalSpace k] [ContinuousAdd k] [ContinuousConstSMul k k] [AddCommMonoid E] [Module k E] [TopologicalSpace E] := WeakBilin (topDualPairing k E)
\end{lstlisting}
\vspace{0.3cm}
Hence, \lean{WeakDual k E} is a type synonym of \lean{E →L[k] k}, equipped with the weak-\(*\) topology.

\begin{figure}[htbp]
    \centering
    \includegraphics[width=0.7\linewidth]{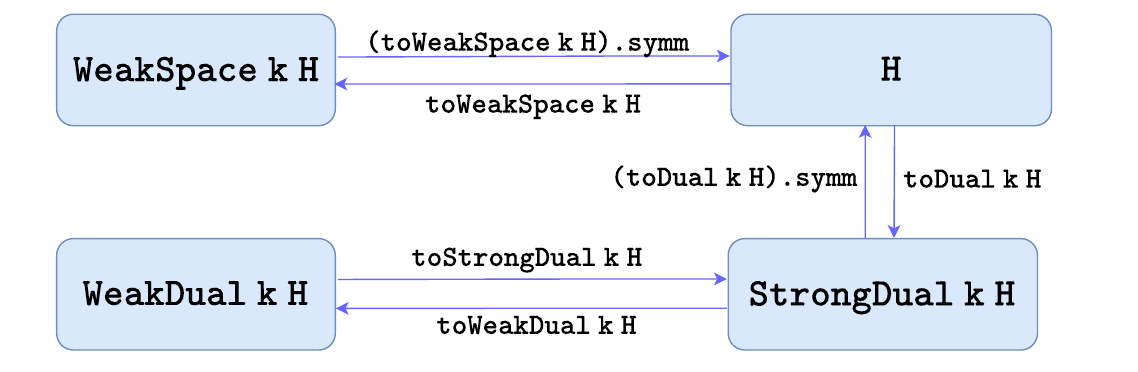}
    \caption{Type coercion between weak and strong topologies of a Hilbert space.}
    \label{fig:weakspace}
\end{figure}

Conversely, the weak topology on $E$ is obtained by using the same bilinear in the opposite order. Concretely, applying \lean{flip} to \lean{topDualPairing k E} yields a pairing of type
\lean{E →ₗ[k] (E →L[k] k) →ₗ[k] k}. The corresponding weak topology on $E$ is formalized as:
\begin{lstlisting}[frame=single]
def WeakSpace (k E) [CommSemiring k] [TopologicalSpace k] [ContinuousAdd k]
    [ContinuousConstSMul k k] [AddCommMonoid E] [Module k E] [TopologicalSpace E] := WeakBilin (topDualPairing k E).flip
\end{lstlisting}

In \mathlib, the dual space $\Hh^*$ of a Hilbert space $\Hh$ over $k$ (real or complex) is denoted as \lean{StrongDual}, which is a type synonym of the continuous linear map \lean{H →L[k] k}. The type transformations among these four spaces in \mathlib are concluded in \cref{fig:weakspace}.

Currently, a wide range of standard results concerning closedness and compactness are formalized for \lean{WeakDual}, such as the Banach–Alaoglu theorem (see \lean{WeakDual.isCompact_closedBall}). However, for Hilbert spaces, \mathlib does not provide a useful coercion between \lean{WeakSpace} and \lean{WeakDual}, which limits the availability of corresponding results in \lean{WeakSpace}. As for weak convergence, although it can be expressed via the inherited, filter-based notion of convergence, convergence in \lean{WeakSpace} is not stated explicitly.

\subsection{Weak convergence}
\label{sec:weak}

As mentioned above, existing topological results for Hilbert spaces in \mathlib mainly focus on the weak-\(*\) topology (\lean{WeakDual}) rather than the weak topology (\lean{WeakSpace}). In this section, we present our formalization in \lean{WeakSpace}. 

Still, let $\Hh$ denotes a real Hilbert space. 
In Lean4, we assume globally that $\Hh$ is a real inner product space, and only add the typeclass \lean{[CompleteSpace H]} to those theorems where completeness is required; in most cases, this dependency arises from invoking the Riesz representation theorem.

\begin{lstlisting}[frame=single]
variable {H : Type*} [NormedAddCommGroup H] [InnerProductSpace ℝ H]
\end{lstlisting}

Most of the proofs in this subsection are based on~\cite[Chapter 2]{bauschke2011convex}. In \mathlib, convergence is defined in the filter theory, hence weak convergence can be expressed by viewing the sequence in the real Hilbert space endowed with the weak topology (\lean{WeakSpace ℝ H}). Recall the type transformation between \lean{H} and \lean{WeakSpace ℝ H} as illustrated in \cref{fig:weakspace}, we map the sequence into \lean{WeakSpace ℝ H} via the composition with \lean{toWeakSpace ℝ H} and then state weak convergence in the weak topology:
\begin{lstlisting}[frame=single]
def WeakConverge (x : ℕ → H) (p : H) :=
  Tendsto ((toWeakSpace ℝ H) ∘ x) atTop (nhds ((toWeakSpace ℝ H) p))
\end{lstlisting}

Throughout, we write $x_n \rightharpoonup p$ in $\Hh$ to denote that $\{x_n \}$ converges weakly to $p$. In many textbooks, weak convergence in a real Hilbert space is characterized in terms of inner product:
\begin{theorem}[Weak convergence]
\label{thm:equiv_weak_conv}
    Let $\{x_n\}$ be a sequence in a real Hilbert space $\Hh$ and let $p \in \Hh$. Sequence $\{x_n\}$ \emph{converges weakly} to $p$ if and only if
    \[
        \ip{x_n}{y} \to \ip{p}{y} \quad \text{for all } y \in \Hh.
    \]
\end{theorem}
We formalize the equivalence between these two definitions as follows, which can be derived by using  \lean{WeakBilin.tendsto_iff_forall_eval_tendsto}. To specialize the continuous linear functionals to the inner product characterization, we use the Riesz representation theorem where the completeness is required.
\begin{lstlisting}[frame=single]
theorem weakConverge_iff_inner_converge [CompleteSpace H] (x : ℕ → H) (p : H) : WeakConverge x p ↔ ∀ y : H, Tendsto (fun n ↦ inner (x n) y) atTop (nhds (inner p y))
\end{lstlisting}

It is well known that strongly convergent sequences are bounded in infinite-dimensional Hilbert spaces. The same property holds for weakly convergent sequences.
\begin{theorem}[Boundedness of weak convergence sequence]
\label{thm:bdd_weak_converge}
If $x_n\rightharpoonup p$, then $\{x_n\}$ is bounded.
\end{theorem}

To prove this claim, for any $n \in \mathbb{N}$, define a bounded linear operator:
\begin{equation*}
     f_n:\Hh \to \mathbb{R}, \quad f_n(y) = \ip{x_n}{y}.
\end{equation*}
Weak convergence of $\{x_n\}$ implies that $\{f_n(y)\}$ is bounded for any $y \in \Hh$. Therefore, the Banach--Steinhaus theorem (in \mathlib, \lean{banach_steinhaus}) yields a uniform bound on $\{f_n\}$. Finally we conclude that $\{ x_n\}$ is bounded since $\norm{x_n} = \norm{f_n}$ in a real Hilbert space. The Lean4 formalization of this theorem is stated as:
\begin{lstlisting}[frame=single]
theorem weakly_converge_norm_bounded [CompleteSpace H] (x : ℕ → H) (p : H)
  (h_wkconv_x : WeakConverge x p) : ∃ M, ∀ n, ‖ x n ‖ ≤ M
\end{lstlisting}

The next theorem shows the weak lower-semicontinuity of the norm.
\begin{theorem}[{\cite[Lemma~2.42]{bauschke2011convex}}]\label{thm:2.42}
If $x_n\rightharpoonup p$, then $\norm{p} \leq {\liminf}_{n\to \infty} \norm{x_n}$.
\end{theorem}
In \mathlib, \lean{liminf} takes values in the extended real numbers \lean{EReal} rather than $\mathbb{R}$. Therefore, we lift $\norm{p}$ from $\mathbb{R}$ to \lean{EReal} using \lean{Real.toEReal}. This theorem is a direct corollary of the Cauchy--Schwarz inequality. The corresponding formal statement is:
\begin{lstlisting}[frame=single]
theorem norm_weakly_lsc [CompleteSpace H] (x : ℕ → H) (p : H) (h : WeakConverge x p) :
  Real.toEReal ‖ p ‖ ≤ liminf (fun n => Real.toEReal ‖ x n ‖) atTop
\end{lstlisting}

The following result is a convenient “mixed convergence” principle: a weakly convergent sequence paired with a strongly convergent sequence yields convergence of the corresponding inner products.
\begin{theorem}[{\cite[Lemma~2.51]{bauschke2011convex}}]\label{thm:2.51}
Suppose $\{x_n\}$ is bounded, and $x_n \rightharpoonup p$, $u_n \to u$, then $\ip{x_n}{u_n}\to\ip{p}{u}$.
\end{theorem}
In \cite{bauschke2011convex}, this theorem is proved using net convergence, which is not convenient in Lean4. Instead, we avoid the net convergence by using the boundedness of $\{x_n\}$ (which follows automatically from weak convergence by \cref{thm:bdd_weak_converge}) and the strong convergence $u_n\to u$. Formalization of this theorem is as follows
\begin{lstlisting}[frame=single]
theorem mix_convergence [CompleteSpace H] {x : ℕ → H} {x_lim : H} {u : ℕ → H} 
{u_lim : H} {h_wkconv_x : WeakConverge x x_lim} {h_conv_u : Tendsto u atTop (nhds u_lim)} : 
  Tendsto (fun n => inner (x n) (u n)) atTop (nhds (inner x_lim u_lim))
\end{lstlisting}

In a finite-dimensional Hilbert space, the weak and strong (norm) topologies coincide. Hence the weak convergence is equivalent to the strong convergence.  The distinction becomes essential only in infinite-dimensional spaces: strong convergence requires \emph{norm} convergence (convergence of $\{\norm{x_n}\}$), whereas weak convergence is strictly weaker.
\begin{theorem}[{\cite[Corollary~2.52]{bauschke2011convex}}]
\label{thm:strong_conv}
In finite-dimensional Hilbert space, the weak convergence is equivalent to the strong convergence. In infinite-dimensional Hilbert space, for any sequence $\{x_n\} \subseteq \Hh$ and $p \in \Hh$, $x_n \to p$ if and only if $x_n \rightharpoonup p$ and $\norm{x_n} \to \norm{p}$.
\end{theorem}
These results are formalized as follows:

\begin{lstlisting}[frame=single]
theorem finite_weak_converge_strong_converge [FiniteDimensional ℝ H] (x : ℕ → H) (p : H) (h : WeakConverge x p) : Tendsto x atTop (nhds p)
theorem strong_converge_weak_converge [CompleteSpace H] (x : ℕ → H) (p : H) (h : Tendsto x atTop (nhds p)) : WeakConverge x p
theorem strong_converge_iff_weak_norm_converge [CompleteSpace H] (x : ℕ → H) (p : H) : Tendsto x atTop (nhds p) ↔ WeakConverge x p ∧ Tendsto (fun n => ‖ x n ‖) atTop (nhds ‖ p ‖)
\end{lstlisting}

\subsection{Closedness}

In \mathlib, topological properties such as \lean{IsClosed} are defined for sets of the ambient type, and the corresponding topological instance is deduced automatically by Lean4. Hence, given a set \lean{s : Set H}, we transform it to the weak topology of $\Hh$ via the canonical map \lean{toWeakSpace ℝ H} and then apply the usual notions of (sequential) closedness there. We define the weak (sequential) closedness in Lean4 as follows:
\begin{lstlisting}[frame=single]
def IsWeaklyClosed (s : Set H) := IsClosed ((toWeakSpace ℝ H) '' s)
def IsWeaklySeqClosed (s : Set H) := IsSeqClosed ((toWeakSpace ℝ H) '' s)
\end{lstlisting}

For convex sets, the weak and strong (sequential) closedness are all equivalent:
\begin{theorem}[{\cite[Theorem~3.34]{bauschke2011convex}}]\label{thm:3.34}
Let $C$ be a convex subset of $\Hh$. The following statement are equivalent:
\begin{enumerate}
    \item $C$ is weakly sequentially closed;
    \item $C$ is sequentially closed;
    \item $C$ is closed;
    \item $C$ is weakly closed.
\end{enumerate}
\end{theorem}
The implication \BLUE{1} $\to$ \BLUE{2} relies on the fact that strong convergence implies weak convergence mentioned in \cref{thm:strong_conv}.
The equivalence \BLUE{2} $\leftrightarrow$ \BLUE{3} is provided by \lean{isSeqClosed_iff_isClosed} in \mathlib.
Convexity is only needed for \BLUE{3} $\to$ \BLUE{4}: it follows from the general theorem \lean{Convex.toWeakSpace_closure} proved in locally convex space, which states that for convex sets, the closure in the strong topology agrees with the closure in the weak topology.
Finally, \BLUE{4} $\to$ \BLUE{1} is immediate from \lean{IsClosed.isSeqClosed}. The equivalence is formalized through \lean{TFAE} (the following are equivalent):
\begin{lstlisting}[frame=single]
theorem seq_closed_tfae [CompleteSpace H] (s : Set H) (hs : Convex ℝ s) :
  [IsWeaklySeqClosed s, IsSeqClosed s, IsClosed s, IsWeaklyClosed s].TFAE
\end{lstlisting}

\subsection{Compactness}
Similar to the definition of closedness, we define compactness in the weak topology of $\Hh$ using the predicate \lean{IsCompact} in \lean{mathlib}:
\begin{lstlisting}[frame=single]
def IsWeaklyCompact (s : Set H) : Prop := IsCompact ((toWeakSpace ℝ H) '' s)
\end{lstlisting}

In general, the dual space of a topological module over a Hausdorff topological commutative semiring is Hausdorff when equiped with the weak-\(*\) topology, which has been formalized as \lean{WeakDual.instT2Space} in \mathlib. However, for a general topological vector space equiped with the weak topology need not be Hausdorff unless the ambient topology is sufficiently separating (e.g., locally convex). In the real Hilbert space setting as considered in this paper, the weak topology is indeed Hausdorff. 

We show in Lean4 that \lean{WeakSpace ℝ H} carries a \lean{T2Space} (Hausdorff) instance:
\begin{lstlisting}[frame=single]
instance inst_WeakSpace_T2 : T2Space (WeakSpace ℝ H)
\end{lstlisting}
Based on this instance, one can obtain the weak closedness of a weak compactness set since in a Hausdorff space, every compact set is closed.

When studying compactness in Hilbert spaces, the Banach--Alaoglu theorem is a pivotal tool in infinite-dimensional analysis and is crucial for many existence and convergence arguments. In particular, this theorem states that every closed ball $B(x, r)$ is weakly compact in a Hilbert space. The weak compactness is fundamental in variational methods and the analysis of weak convergence.
\begin{theorem}[Banach--Alaoglu theorem, {\cite[Fact~2.34]{bauschke2011convex}}]
\label{thm:Banach}
The closed ball $B(x,r)$ is weakly compact.
\end{theorem}

In \mathlib, a similar theorem is formalized in \lean{WeakDual}: the closed ball in $\Hh^*$ is weakly compact, see \lean{WeakDual.isCompact_closedBall}. We formalize the corresponding version in \lean{WeakSpace ℝ H} as follows:
\begin{lstlisting}[frame=single]
theorem closed_unit_ball_is_weakly_compact [CompleteSpace H] (x : H) (r : ℝ) :
  IsWeaklyCompact (closedBall x r)
\end{lstlisting}

According to the Riesz representation theorem, $\Hh$ is not that different from $\Hh^*$. The key step in relating their weak topologies is to establish a type conversion between the two spaces. However, such a coercion is not currently available in \mathlib, see \cref{fig:weakspace}. To accomplish this, we construct an explicit homeomorphism between \lean{WeakSpace ℝ H} and \lean{WeakDual ℝ H}, as shown in the following Lean4 code:
\begin{lstlisting}[frame=single]
noncomputable def weakHomeomorph [CompleteSpace H] : Homeomorph (WeakSpace ℝ H) (WeakDual ℝ H)
\end{lstlisting}

Since homeomorphism preserves compactness (\lean{Homeomorph.isCompact_image} in \mathlib), we complete the proof for \cref{thm:Banach}.

To ensure the existence of a weakly convergent subsequence, it suffices to show the boundedness.
\begin{theorem}[{\cite[Lemma~2.45]{bauschke2011convex}}]
\label{Lemma 2.45}
Every bounded sequence $\{x_n\}$ in $\Hh$ possesses a weak subsequential limit point, which is the limit point of a weakly convergent subsequence.
\end{theorem}

The weak subsequential limit point of a sequence is formalized as below. The subsequence is expressed through a strict monotone index sequence.
\begin{lstlisting}[frame=single]
def WeakSubseqLimitPt (p : H) (x : ℕ → H) :=
 ∃ (φ : ℕ → ℕ), StrictMono φ ∧ WeakConverge (x ∘ φ) p
\end{lstlisting}

To facilitate the Lean4 formalization of \cref{Lemma 2.45}, we prove a 
specialized version under the assumption that the Hilbert space $\Hh$ is separable. We follow the diagonal argument method in~\cite{royden1988real} by explicitly constructing a weak convergent subsequence. The Lean4 statement corresponding to \cref{Lemma 2.45} is:
\begin{lstlisting}[frame=single]
theorem bounded_seq_has_WeakSubseqLimitPt_separable [SeparableSpace H]
  [CompleteSpace H] (x : ℕ → H) (hx : Bornology.IsBounded <| Set.range (fun n => ‖ x n ‖)) : ∃ (a : H), WeakSubseqLimitPt a x
\end{lstlisting}

As a corollary, if the set of weak subsequential limit points of $\{x_n\}$ contains at most one element, then it possesses exactly one weak subsequential limit point. Consequently, the whole sequence converges weakly.
\begin{corollary}[{\cite[Lemma~2.47]{bauschke2011convex}}]
\label{thm:2.47}
If sequence $\{x_n\}$ is bounded and possesses at most one weak subsequential limit point, then $\{x_n\}$ weakly converges to some point $p_0\in\Hh$.
\end{corollary}
The formalization of this theorem is as follows:
\begin{lstlisting}[frame=single]
lemma WeakConv_of_bounded_clusterptUnique [SeparableSpace H] [CompleteSpace H] (x : ℕ → H) (hb : ∃ M : ℝ, ∀ n, ‖ x n ‖ ≤ M)
(hc : ∀ p q : H,  WeakSubseqLimitPt p x → WeakSubseqLimitPt q x  → p = q) : ∃ p0 : H, WeakConverge x p0
\end{lstlisting}

\section{Nonexpansive operator and fixed-point iteration}
\label{sec:nonexp}

Nonexpansive operators are ${1}$-Lipschitz mappings, meaning that such operators will not increase distance. In computational mathematics, many fundamental iterative methods can be written as repeated applications of a (or closely related) nonexpansive operator. Nonexpansiveness provides convenient geometric tools for stability and convergence analysis, and it also underlies convergence results for operator splitting methods such as Douglas–Rachford splitting~\cite{lions1979splitting}, and alternating direction method of multipliers (ADMM)~\cite{glowinski1975approximation, gabay1983chapter}, which can be analyzed within the framework of nonexpansive operators~\cite{eckstein1992douglas}.

In this section, we formalize nonexpansive operators and the two fixed-point iteration schemes studied in this paper, each defined as an iterate of a nonexpansive operator.

\subsection{Nonexpansive operators}

\begin{definition}[Nonexpansive operator]
\label{def:nonexp}
Let $D \subseteq \Hh$, and let $T : \Hh \to \Hh$ be an operator. We say that $T$ is \emph{nonexpansive on $D$} if, for all $x,y \in D$,
\[
\|Tx - Ty\| \le \|x - y\|.
\]
If $D=\Hh$, we simply say that $T$ is \emph{nonexpansive}.
\end{definition}
As a special case of Lipschitz continuity, we use \lean{LipschitzWith} and \lean{LipschitzOnWith} in \mathlib to define nonexpansiveness on the pseudo extended metric space:
\begin{lstlisting}[frame=single]
variable [PseudoEMetricSpace α] [PseudoEMetricSpace β]
def Nonexpansive (T : α → β) := LipschitzWith 1 T
def NonexpansiveOn (T : α → β) (s : Set α) := LipschitzOnWith 1 T s
\end{lstlisting}

In certain scenarios, one may only concerns the distance to a fixed-point being controlled. Consequently, the nonexpansiveness reduces to  \emph{quasinonexpansiveness}.

\begin{definition}[Quasinonexpansive]
Let $D \subseteq \Hh$, and let $T : \Hh \to \Hh$ be an operator. We say that $T$ is \emph{quasinonexpansive on $D$} if, for all $x \in D, y \in \Fix T\cap D$,
\[
\|Tx - y\| \le \|x - y\|.
\]
If $D=\Hh$, we simply say that $T$ is \emph{quasinonexpansive}.
\end{definition}

We first define the fixed-point set of an operator $T$ as follows by using \lean{IsFixedPt}:
\begin{lstlisting}[frame=single]
def Fix (T : α → α) : Set α := {x | IsFixedPt T x}
\end{lstlisting}

In stead of restricting to Hilbert space, we work in the more general normed additive commutative group (a space equipped with a norm and an additive commutative group structure).
\begin{lstlisting}[frame=single]
variable [NormedAddCommGroup α]
def QuasiNonexpansive (T : α → α) := ∀ x y, y ∈ Fix T → ‖ T x - y ‖ ≤ ‖ x - y ‖
def QuasiNonexpansiveOn (T : α → α) (s : Set α) :=
  ∀ ⦃x⦄, x ∈ s → ∀ ⦃y⦄, y ∈ Fix T ∩ s → ‖ T x - y ‖ ≤ ‖ x - y ‖
\end{lstlisting}

Clearly, every nonexpansive operator is quasinonexpansive, formalized as follows:
\begin{lstlisting}[frame=single]
theorem nonexpansive_quasinonexpansive {D : Set α} {T : α → α}
  (hT_nonexp : NonexpansiveOn T D) : QuasiNonexpansiveOn T D
\end{lstlisting}

In Hilbert spaces, it can be shown that when $D$ is closed and convex, the fixed-point set of a quasinonexpansive operator is closed and convex.
\begin{theorem}[{\cite[Proposition~4.23]{bauschke2011convex}}]
Suppose $D$ is a nonempty closed and convex subset of $\Hh$. If $T:\Hh \to \Hh$ is quasinonexpansive on $D$, then $\Fix T \cap D$ is closed and convex.
\end{theorem}
To prove this result, we first provide a characterization of the fixed-point set as an intersection of closed and convex half-spaces, from which the claim follows immediately.
\begin{equation*}
    \Fix T \cap D = \bigcap_{x \in D} \Big\{y \in D \mid \ip{y - T x}{x - T x} \leq \frac{1}{2} \norm{T x - x}^2\Big\}.
\end{equation*}
The above two theorems are formalized in Lean4 as follows:
\begin{lstlisting}[frame=single]
theorem quasinonexpansive_fixedPoint_characterization {D : Set H} {T : H → H} (hD_nonempty : D.Nonempty) (hT_quasi : QuasiNonexpansiveOn T D) :
  Fix T ∩ D = ∩ x ∈ D, {y ∈ D | inner (y - T x) (x - T x) ≤ (1/2) * ‖ T x - x ‖^2} 

theorem quasinonexpansive_fixedPoint_closed_convex {C D : Set H} (hD_closed : IsClosed D) (hD_convex : Convex ℝ D) (hD_nonempty : D.Nonempty) {T : H → H} 
(hT : QuasiNonexpansiveOn T D) (hC : C = Fix T ∩ D) : IsClosed C ∧ Convex ℝ C
\end{lstlisting}

Next we introduce another special operator called demiclosed operator, which describes the image convergence of an operator over a weakly convergent sequence.

\begin{definition}[Demiclosed operator]
\label{def:demiclosed}
Let $D$ be a nonempty weakly sequentially closed subset of $\Hh$, and let $u \in \Hh$. Then $T:\Hh \to \Hh$ is demiclosed at $u$ if for every sequence $\{x_n \}$ in $D$ such that $x_n \rightharpoonup x $ and $Tx_n \to u$, we have $Tx = u$. In addition, $T$ is demiclosed on $D$ if it is demiclosed at every point in $D$.
\end{definition}

\begin{lstlisting}[frame=single]
def DemiclosedAt (D : Set H) (T : H → H) (u : H) : Prop :=
  (hD : D.Nonempty) → (hw : IsWeaklySeqClosed D) →
  ∀ (x : ℕ → H), (∀ n, x n ∈ D) → ∀ (l : H), l ∈ D → WeakConverge x l →
  Tendsto (fun n => T (x n)) atTop (nhds u) → T l = u

def DemiclosedOn (T : H → H) (D : Set H) := ∀ u : H, DemiclosedAt D T u
\end{lstlisting}

A fundamental fact known as Browder's demiclosedness principle is that for a nonexpansive operator $T$, $\Id - T$ is demiclosed.

\begin{theorem}[{\cite[Theorem~4.27]{bauschke2011convex}}]
\label{thm:demiclosed}
Let $D$ be a nonempty weakly sequentially closed subset of $\Hh$, and $T:\Hh \to\Hh$ be a nonexpansive operator on $D$. Then $\Id - T$ is demiclosed.
\end{theorem}

\begin{lstlisting}[frame=single]
theorem browder_demiclosed_principle [CompleteSpace H] {D : Set H} {T : H → H} (hT_nonexp : NonexpansiveOn T D) : DemiclosedOn (id - T) D
\end{lstlisting}

As a direct consequence, since
\begin{equation*}
    (\Id-T)x=0 \iff x \in \mathrm{Fix}\,T,
\end{equation*}
\cref{thm:demiclosed} provides a convenient criterion for verifying that a weak limit point belongs to $\Fix T$, which is concluded in the following corollary:

\begin{corollary}[{\cite[Corollary~4.28]{bauschke2011convex}}]
\label{lem:x-in-FixT}
Let $D$ be a nonempty closed and convex subset of $\Hh$, let $T : \Hh \to \Hh$ be nonexpansive on $D$, let $\{x_n\}$ be a sequence in $D$, and let $x\in D$. Suppose that $x_n \rightharpoonup x$ and that $x_n - Tx_n \to 0$. Then $x \in \Fix T$.
\end{corollary}
The formalization of this theorem is as follows:
\begin{lstlisting}[frame=single]
lemma weakLimit_mem_fixedPoints_of_strongly_tendsto_sub [CompleteSpace H]
  {D : Set H} (hD_closed : IsClosed D) (hDc : Convex ℝ D) (hDn : D.Nonempty) {T : H → H} 
  (hT_nonexp : NonexpansiveOn T D) (x : ℕ → H) (hx : ∀ n, x n ∈ D) (p : H) (hp : p ∈ D) 
  (hw : WeakConverge x p) 
  (he : Tendsto (fun n => x n - T (x n)) atTop (nhds 0)) : p ∈ Fix T
\end{lstlisting}

\subsection{Fixed-point iterations}
\label{sec:fixed-point}

In this part, suppose $D \subseteq \mathcal{H}$ is a nonempty closed convex set, and $T: \Hh \to \Hh$ is a nonexpansive operator with $T(D) \subseteq D$. Assume that the fixed-point set $C \triangleq \Fix T \cap D$ is nonempty.

Recall the update scheme of KM iteration \cref{eq:KM},
\begin{equation}
\forall n \in \mathbb{N}, \quad x_{n+1} = x_n + \alpha_n(Tx_n - x_n),
\end{equation}
where the initial point $x_0$ belongs to $D$. We use Lean4's \lean{structure} to define the iterative scheme of the KM iteration as follows.

\begin{lstlisting}[frame=single]
structure KM (D : Set H) (T : H → H) where
  x0 : H   -- initial point
  α : ℕ → ℝ   -- relaxation parameters
  x : ℕ → H   -- sequence of iterates
  hx0 : x0 ∈ D      
  update : ∀ n, x (n + 1) = x n + α n • (T (x n) - x n)
  initial_value : x 0 = x0
\end{lstlisting}
The \lean{KM} structure bundles the domain \(D\), operator \(T\), iterate sequence \(\{x_n\}\), initial value \(x_0\). Thus, a term \lean{km : KM D T} constructs a concrete KM iteration instance, where \lean{km.x} denotes the iterate sequence and \lean{km.update} records the update scheme. This structured approach improves the modularity and reusability of formalization.

We then formalize the Halpern iteration scheme \cref{eq:Halpern}, restated here for convenience:
\begin{equation*}
    \forall n \in \mathbb{N}, \quad x_{n+1} = \alpha_n u + (1 - \alpha_n) T x_n.
\end{equation*}
Here $u,x_0 \in D$. The formalization of the above Halpern iteration is as follows, using \lean{structure} to package all the components:
\begin{lstlisting}[frame=single]
structure Halpern (T : H → H) where
  x0 : H   -- initial point  
  u : H   -- anchor point
  x : ℕ → H   -- sequence of iterates
  α : ℕ → ℝ   -- sequence of weight
  update : ∀ k : ℕ, x (k + 1) = (α k) • u + (1 - α k) • (T (x k))
  initial_value : x 0 = x0
\end{lstlisting}

\section{Convergence analysis}
\label{sec:convergence}

Leveraging the formalized scheme in \cref{sec:fixed-point}, we introduce the main convergence results in this section.

\subsection{Convergence of \KM iteration}
\label{sec:KM_conv}

We first present a formalization of weak convergence of the classical \KM (KM) iteration algorithm in Lean4. Our development builds on the notions of weak convergence and nonexpansive operators introduced in \cref{sec:Weak,sec:nonexp}. To establish weak convergence of the iterates, we introduce \Fejer monotonicity as a key tool and then prove the main convergence theorem in \cref{sec:KM}.

\subsubsection{\Fejer monotonicity}
\label{sec:fejer}

Fej\'{e}r monotonicity is a widely used tool in fixed-point theory. Generally speaking, it describes a property of a sequence $\{x_n \}$ that every iterate $x_{n+1}$ is no farther from any point in a subset $D \subseteq \Hh$ than the previous iterate $x_{n}$. As what we will present later, Fej\'{e}r monotonicity immediately implies boundedness of $\{x_n\}$ and provides a strong geometric control on the iterates. Combined with mild additional conditions (e.g., weak sequential compactness and suitable asymptotic regularity), it often leads to convergence of the sequence to a point in $D$.

\begin{definition}[Fej\'{e}r monotone]
\label{def:fejer}
A sequence $\{x_n \}$ in $\Hh$ is said to be Fej\'{e}r monotone with respect to a set $D \subseteq \Hh$ if for any $x \in D$,
\begin{equation*}
\forall n \in \mathbb{N}, \qquad \norm{x_{n+1} - x } \leq \norm{x_n - x}.
\end{equation*}
\end{definition}
The formalization of this definition is as follows:
\begin{lstlisting}[frame=single]
def IsFejerMonotone (x : ℕ → H) (D : Set H) : Prop :=
  ∀ y ∈ D, ∀ n, ‖ x (n + 1) - y ‖ ≤ ‖ x n - y ‖
\end{lstlisting}

A direct observation is that Fej\'{e}r monotone sequences are bounded. Moreover, for any $y \in D$, the sequence $\{\norm{x_n-y}\}$ is convergent since it is monotonically decreasing and bounded. We formalize these two basic results as follows:

\begin{lstlisting}[frame=single]
theorem Fejermono_bounded (D : Set H) (hD : D.Nonempty) (x : ℕ → H)
  (hx : IsFejerMonotone x D) : ∃ M:ℝ , ∀ n, ‖ x n ‖ ≤ M
  
theorem Fejermono_convergent (D : Set H) (x : ℕ → H) (h : IsFejerMonotone x D) : ∀ a ∈ D, ∃ l : ℝ, Tendsto (fun n ↦ ‖ x n - a ‖) atTop (nhds l)
\end{lstlisting}

Based on the convergence of $\{\norm{x_n -a} \}$, the weak convergence of Fej\'{e}r monotone sequence can be derived in the following theorem:
\begin{theorem}
\label{thm:5.5}
If the sequence $\{x_n\}$ is Fej\'{e}r monotone with respect to a nonempty set $D$, and if all the weak sequential limit points of $\{x_n \}$ belong to $D$, then
$\{x_n\}$ weakly converges to some point $p_0$ in $D$.
\end{theorem}
By \cref{thm:2.47}, to prove the weak convergence of $\{x_n\}$, it suffices to show that $\{x_n\}$ has exactly one weak sequential limit point, which can be deduced from the Fej\'{e}r monotonicity of the sequence $\{x_n\}$. The corresponding Lean4 statement of \cref{thm:5.5} is provided below.

\begin{lstlisting}[frame=single]
theorem WeakConv_of_Fejermonotone_of_clusterpt_in [SeparableSpace H] [CompleteSpace H] (D : Set H) (hD : D.Nonempty) (x : ℕ → H)
(hx : IsFejerMonotone x D) (hw : ∀ p : H, WeakSubseqLimitPt p x → p ∈ D):
  ∃ p0 ∈ D, WeakConverge x p0
\end{lstlisting}

\subsubsection{Convergence}
\label{sec:KM}

To guarantee the weak convergence of the generated sequence $\{x_n\}$, the relaxation parameter $\{\alpha_n \}$ is assumed to satisfy the following classical conditions.
\begin{assumption}\label{ass:KM}
The relaxation parameter $\alpha_n$ satisfies
\begin{enumerate}
	\item $\alpha_n \in [0,1]$.
	\item $\displaystyle\sum_{n=0}^{\infty} \alpha_n(1 - \alpha_n) = +\infty$.
\end{enumerate}
\end{assumption}
Condition \BLUE{1} ensures that the iterates remain in the domain $D$ by convexity. Condition \BLUE{2} guarantees that the iteration does not become asymptotically stagnant.

The convergence of the KM iteration is established by the following theorem. Since Groetsch~\cite{GROETSCH1972369} provided sufficient conditions for the convergence of the KM algorithm, we denote our main theorem as \lean{Groetsch_theorem}.

\begin{theorem}[Convergence of KM Iteration, {\cite[Theorem~5.15]{bauschke2011convex}}]
\label{thm:groetsch}
Denote $\{x_n \}$ the sequence generated by the KM iteration \eqref{eq:KM}, and 
suppose \cref{ass:KM} holds. Let $C \triangleq \Fix T \cap D$.
Then the following claims hold:
	\begin{enumerate}
		\item $\{x_n\}_{n \in \mathbb{N}}$ is Fej\'er monotone with respect to $C$.
		\item $\{Tx_n - x_n\}_{n \in \mathbb{N}}$ converges strongly to $0$.
		\item $\{x_n\}_{n \in \mathbb{N}}$ converges weakly to a point in $C$.
	\end{enumerate}
\end{theorem}
\vspace{0.9cm}
The formalized version of the theorem is:
\begin{lstlisting}[frame=single]
theorem Groetsch_theorem [SeparableSpace H] [CompleteSpace H] {D : Set H}
  (hD1 : Convex ℝ D) (hD2 : IsClosed D) (T : H → H) (h_Im_T_in_D : ∀ x ∈ D, T x ∈ D)
  (hT : NonexpansiveOn T D) (km : KM D T) (fix_T_nonempty : (Fix T ∩ D).Nonempty)
  (hα1 : ∀ n, km.α n ∈ Set.Icc (0 : ℝ) 1)
  (hα2 : Tendsto (fun n => ∑ i ∈ range (n + 1), km.α i * (1 - km.α i)) atTop atTop) :
    IsFejerMonotone km.x (Fix T ∩ D)
    ∧ (Tendsto (fun n ↦ ‖ T (km.x n) - km.x n ‖) atTop (nhds 0))
    ∧ ∃ y0 ∈ (Fix T ∩ D), WeakConverge km.x y0
\end{lstlisting}

We mainly follow the proof of~\cite[Theorem 5.15]{bauschke2011convex}. Since the hypotheses that $D$ is closed and convex and that $T$ is nonexpansive are not part of the algorithm but rather analytical prerequisites, they are explicitly stated in the theorem and labeled with descriptive identifiers such as \lean{hD_closed} and \lean{hT}.
The proof proceeds in three stages: 
\begin{enumerate}
    \item The Fej\'{e}r monotonicity of  sequence $\{x_n\}$ with respect to $C$ can be derived by the following inequality:
\begin{equation*}
    	\forall y \in C, \qquad \|x_{n+1} - y\|^2 
	\leq \|x_n - y\|^2 - \alpha_n(1 - \alpha_n)\|Tx_n - x_n\|^2,
\end{equation*}
which is formalized as below:
\begin{lstlisting}[frame=single]
lemma key_inequality {D : Set H} (T : H → H) (h_Im_T_in_D : ∀ x ∈ D, T x ∈ D)
  (hT : NonexpansiveOn T D) (km : KM D T) (hD : Convex ℝ D) 
  (hα1 : ∀ n, km.α n ∈ Set.Icc (0 : ℝ) 1)
  (hα2 : Tendsto (fun n => ∑ i ∈ range (n + 1), km.α i * (1 - km.α i)) atTop atTop) :
    ∀ (y : H) (hy : y ∈ Fix T ∩ D) (n : ℕ), ‖ km.x (n + 1) - y ‖^2 ≤ ‖ km.x n - y ‖^2 - km.α n * (1 - km.α n) * ‖ T (km.x n) - km.x n ‖^2
\end{lstlisting} 
 
\item Furthermore, based on the inequality above and \cref{ass:KM} \BLUE{2}, one can derive that
\begin{equation*}
    \liminf \norm{Tx_n - x_n} = 0.
\end{equation*}
On the other hand, we prove that $\{\norm{Tx_n - x_n}\}$ is monotonically decreasing, hence $\norm{Tx_n - x_n} \to 0$.
\item To establish weak convergence of the Fej\'{e}r monotone sequence $\{x_n\}$, it suffices to prove that all the weak subsequential limit points of $\{ x_n\}$ belongs to $C$ according to \cref{thm:5.5}.
Let $x$ be be such a weak subsequential limit point of $\{x_n\}$, say $x_{k_n} \rightharpoonup x$. Then it follows from \cref{lem:x-in-FixT} and \cref{thm:groetsch} \BLUE{2} that $x \in C$. Therefore, $\{x_n \}$ converges weakly to a point in $C$.

\end{enumerate}

This completes the formal proof of convergence for the KM algorithm. Among the convergence assumptions, the closedness of \(D\) is required only in the 3rd claim of \cref{thm:groetsch}, while the separability and completeness of the Hilbert space \(\Hh\) are invoked solely in the final step to guarantee the existence of weak limit point according to \cref{thm:5.5}.

\subsection{Convergence of Halpern iteration}
\label{sec:halpern}

To guarantee the convergence of Halpern iteration, the weight $\alpha_n$ between the anchor point $u$ and the image $Tx_n$ should satisfy the following assumptions:
\begin{assumption}\label{ass:Halpern}
The relaxation parameter $\alpha_n$ satisfies
\begin{enumerate}
    \item For any $ n \in \mathbb{N},\alpha_{n}\in(0,1)$;
    \item $\{ \alpha_n\}$ converges to $0$;
    \item $\sum_{n\in\mathbb{N}}\alpha_{n}=+\infty$;
    \item $\sum_{n\in\mathbb{N}}|\alpha_{n+1}-\alpha_{n}|<+\infty$. 
\end{enumerate}
\end{assumption}
Condition \BLUE{1} reflects the convex combination nature of Halpern iteration. Condition \BLUE{2} ensures that the influence of the anchor point vanishes asymptotically. Condition \BLUE{3} guarantees that the cumulative influence of the anchor point is strong enough. Condition \BLUE{4} controls the oscillation of the parameters.
Under these assumptions, the sequence generated by Halpern iteration is guaranteed to converge to a fixed-point of $T$.

\begin{theorem}[Convergence of Halpern Iteration, {\cite[Theorem~30.1]{bauschke2011convex}}]
\label{thm:halpern}
Denote $\{x_n \}$ the iterate sequence generated by Halpern iteration \eqref{eq:Halpern}. Let $C \triangleq \Fix T \cap D$. Assume that \cref{ass:Halpern} holds, and $C$ is nonempty, then $x_n \to P_{C} u$, where $P_{C}u$ denotes the projection point of $u$ to $C$.
\end{theorem}
The formalization of this main theorem is as follows:
\begin{lstlisting}[frame=single]
theorem halpern_convergence [CompleteSpace H] [SeparableSpace H]
  {D : Set H} (hD_closed : IsClosed D) (hD_convex : Convex ℝ D) (hD_nonempty : D.Nonempty) {T : H → H} (hT_nonexp : NonexpansiveOn T D) 
  {C : Set H} (hC : C = Fix T ∩ D) (hCn : C.Nonempty)
  (hT_inD : ∀ x ∈ D, T x ∈ D) (alg : Halpern T)
  (halg_x0 : alg.x0 ∈ D) (halg_u : alg.u ∈ D) (halg_xD : ∀ n, alg.x n ∈ D)
  (hα1 : ∀ n, alg.α n ∈ Set.Ioo 0 1) 
  (hα2 : Tendsto alg.α atTop (nhds 0))
  (hα3 : Tendsto (fun N => ∑ n ∈ Finset.range N, alg.α n) atTop atTop)
  (hα4 : Summable (fun n => |alg.α (n + 1) - alg.α n|)) :
  ∃ (p : H), p ∈ C ∧ Tendsto alg.x atTop (nhds p) ∧ ( ‖ alg.u - p ‖ = ⨅ w : C, ‖ alg.u - w ‖ )
\end{lstlisting}
To characterize the projection onto a closed and convex set, we use the following two theorems from \texttt{mathlib}:
\begin{itemize}
    \item \lean{exists_norm_eq_iInf_of_complete_convex};
    \item \lean{norm_eq_iInf_iff_real_inner_le_zero}.
\end{itemize}
Accordingly, we write $p=P_{C }u$ in Lean4 by choosing $p \in C$ such that $p$ attains the minimum distance, namely
\begin{equation*}
    \norm{u - p} = \inf_{w \in P_{C}} \norm{u-w}.
\end{equation*}
We mainly follow the proof of~\cite[Theorem 30.1]{bauschke2011convex}, which can be divided into two different cases depending on whether the initial point coincides with the anchor point.

\textbf{Case 1: $x_0=u$.}  The proof is summarized as below.
\begin{enumerate}
    \item We first establish the boundedness of the sequences $\{x_n\}, \{Tx_n\}, \{x_{n+1} - Tx_n\}$ by induction. This step serves two purposes: first, it provides an essential uniform bound required in the subsequent estimates; second, it guarantees that $\{x_n\}$ admits a weakly convergent subsequence according to \cref{Lemma 2.45}, which is a key ingredient in the analysis of the limit inferior.
    
    \item We show that the difference between two iterates is close enough: 
    \begin{equation*}
        \limsup \norm{x_{n+2}-x_{n+1}} \leq 0,
    \end{equation*}
    and hence $x_{n+1}-x_n\to 0$.
    Moreover, by the nonexpansiveness of $T$ we have $Tx_{n+1}-Tx_n\to 0$.
    Combining this with the boundedness of $\{x_{n+1}-Tx_n\}$ yields $x_n-Tx_n\to 0$.
    
    \item 
    To establish the desired convergence, we prove the following angle condition between $\{x_n\}$ and the projection point $P_{C}u$:
    \begin{equation*}
        \limsup \ip{Tx_n - P_{C} u}{u -  P_{C} u} \leq 0.
    \end{equation*}
    The analysis of this limit superior proceeds by extracting a weakly convergent subsequence $\{x_{k_n} \}$ of $\{x_n\}$, whose existence is guaranteed by the boundedness of $\{x_n\}$ and \cref{Lemma 2.45}. This is the main reason why separability of the Hilbert space is assumed in the Lean4 formalization of \cref{thm:halpern},
    since \cref{Lemma 2.45} requires separability.
    
    \item We derive a key estimate showing that the iterates become arbitrarily
    close to the projection point $P_{C}u$ when the number of iterations is large enough: for any $\varepsilon>0$, there
    exists $m\in\mathbb{N}$ such that for all $n\ge m$,
    \begin{equation*}
        \norm{x_{n+1}-P_{C}u}^{2}\leq 3\varepsilon+\norm{x_{m}-P_{C}u }^{2} \cdot 
        \mathrm{exp} \left( -\sum_{k=m}^n \alpha_k\right).
    \end{equation*}
    Based on \cref{ass:Halpern} \BLUE{2}, as $n \to \infty$,
    \begin{equation*}
        \mathrm{exp} \left( -\sum\nolimits_{k=m}^n \alpha_k\right) \to 0,
    \end{equation*}
    therefore we derive the convergence of $\{x_n\}$, that $
        x_{n+1} \to P_{C} u$. 
\end{enumerate}

\textbf{Case 2: $x_0\neq u$.}  In this case, we denote $\{y_n \}$ the sequence generated by Halpern iteration \eqref{eq:Halpern} with initial point $y_0 = u$. By the previous discussion, when the initial point equals to the anchor point, we can deduce that $y_n \to P_{C} u$. For $\{x_n\}$, we prove the following estimate between the two different iterates:
\begin{equation*}
    \forall n \in \mathbb{N}, \qquad \norm{x_{n+1} - y_{n+1}} \leq \norm{x_0 - y_0} \cdot \mathrm{exp} \left( \sum_{k=0}^n -\alpha_k \right) \to 0.
\end{equation*}
Hence $x_n \to P_{C} u$. The inequality above is formalized in \lean{halpern_norm_diff_bound}, and the convergence is formalized in \lean{halpern_prod_norm_diff_tendsto_zero}.
Combining the above two cases, we complete the convergence analysis of Halpern iteration.

\section{Conclusion}
\label{sec:con}
In this paper, we formalize the convergence analysis of two fixed-point iteration schemes for nonexpansive operators: the Halpern iteration and the Krasnosel'ski\u{\i}–Mann (KM) iteration. Accordingly, we formalize the definition and key properties of nonexpansive operators, which play an essential role in convergence proofs.
As a further tool for the analysis, we also formalize Fej\'{e}r monotone sequences and several fundamental properties. Since the KM iteration is commonly studied in infinite-dimensional Hilbert spaces under weak convergence, we additionally develop a collection of basic results concerning the weak topology.
For future work, we plan to formalize convergence rates~\cite{liang2016convergence, lieder2021convergence} for these two iterations and to extend the development to broader classes of operators, such as monotone operators arising in convex analysis. Moreover, while the present formalization is carried out in real Hilbert spaces, we expect that results mentioned in \cref{sec:weak} can be generalized to more general settings.



\section*{Acknowledgments} The authors would like to thank the \lean{Zulip} community for valuable discussion on this work.

\bibliography{reference}

\appendix

\end{document}